\documentclass[11pt,oneside]{amsart}
\usepackage{amsfonts, amscd, amsmath, amssymb}
\usepackage[dvips]{graphicx}

\newcommand\globcnt{subsubsection}

\swapnumbers
\theoremstyle{plain}
\newtheorem{theorem}[\globcnt]{Theorem}         
\newtheorem{proposition}[\globcnt]{Proposition} 
\newtheorem{corollary}[\globcnt]{Corollary}
\newtheorem{lemma}[\globcnt]{Lemma}
\theoremstyle{definition}

\newtheorem{definition}[\globcnt]{Definition}

\newtheorem{remark}[\globcnt]{Remark}

\newtheorem{example}[\globcnt]{Example}

\makeatletter
\@addtoreset{subsubsection}{section}
\@addtoreset{subsubsection}{subsection}
\@addtoreset{equation}{subsection}
\@addtoreset{figure}{subsection}
\makeatother

\newcommand\chEq{%
\setcounter{equation}{\value{subsubsection}}%
\addtocounter{subsubsection}{1}
}

\newcommand\chFig{%
\setcounter{figure}{\value{subsubsection}}%
\addtocounter{subsubsection}{1}
}

\title {Stabilizers and orbits of circle-valued smooth functions} 
\author{Sergey Maksymenko}

\address{Topology Dept., Institute of Mathematics, NAS of Ukraine,   Tere\-shchenkivs\-ka str. 3, 01601 Kyiv, Ukraine, e-mail:\texttt{maks@imath.kiev.ua}}

\providecommand{\eqref}[1]{(\ref{#1})}

\newcommand\CCC{{\mathbb C}}

\newcommand\RRR{{\mathbb R}}

\newcommand\ZZZ{{\mathbb Z}}

\newcommand\FFF{{\mathcal F}}

\providecommand\proof{\proofstyle{Proof.\ }}

\providecommand{\id}{\mathop{\rm id}\nolimits}

\newcommand\Orbit{\mathcal{O}}
\newcommand\Stab{\mathcal{S}}
\newcommand\Diff{\mathcal{D}}

\newcommand\pnt{z}

\newcommand\eps{\varepsilon}

\newcommand\manif{M}
\newcommand\flow{\Phi}
\newcommand\aflow{\Psi}

\newcommand\difR{\phi}
\newcommand\difP{\difR}

\newcommand\nbh{U}

\newcommand\dimM{m}
\newcommand\difM{h}
\newcommand\Fld{F}

\newcommand\mrsfunc{f}

\newcommand\manpsp{\manif\Psp}
\newcommand\manr{\manif\RRR}
\newcommand\mans{\manif\Circle}

\newcommand\brf{(\mrsfunc)}

\newcommand\DiffMP{\Diff_{\manif\Psp}}
\newcommand\DiffMR{\Diff_{\manif\RRR}}
\newcommand\DiffMS{\Diff_{\manif\Circle}}
\newcommand\DiffP{\Diff_{\Psp}}
\newcommand\DiffR{\Diff_{\RRR}}
\newcommand\DiffS{\Diff_{\Circle}}
\newcommand\DiffM{\Diff_{\manif}}
\newcommand\StabfMP{\Stab_{\manpsp}}
\newcommand\StabfMR{\Stab_{\manr}}
\newcommand\StabfM{\Stab_{\manif}}
\newcommand\StabfMS{\Stab_{\mans}}
\newcommand\pStabfMS{\widetilde{\Stab}_{\mans}}

\newcommand\OrbfMP{\Orbit_{\manpsp}}
\newcommand\OrbfMS{\Orbit_{\mans}}
\newcommand\pOrbfMS{\widetilde{\Orbit}_{\mans}}
\newcommand\OrbfM{\Orbit_{\manif}}
\newcommand\OrbfMR{\Orbit_{\manr}}

\newcommand\DiffSpaceFunc[1]{\Diff_{#1}}
\newcommand\DiffSpaceCrVal[1]{\Diff_{#1}^{e}}
\newcommand\DiffSpaceCrValInv[1]{\Diff_{#1}^{E}}

\newcommand\DiffRf{\DiffSpaceFunc{\RRR}}
\newcommand\DiffSf{\Diff^{+}_{\Circle}}

\newcommand\DiffSCrValf{\DiffSpaceCrVal{\Circle}}
\newcommand\DiffRCrValf{\DiffSpaceCrVal{\RRR}}

\newcommand\DiffSCrValfInv{\DiffSpaceCrValInv{\Circle}}

\newcommand\claimbmap{{\mathbf {sh}}}

\newcommand\Shift{\varphi}
\newcommand\difShift{\Omega}

\newcommand\Psp{P} 
\newcommand\Circle{S^1}
\newcommand\afunc{\alpha}

\newcommand\clmdif{\sigma} 

\newcommand\prMPP{p}

\newcommand\sct{s}

\newcommand\symplex{\Delta^{n-1}}
\newcommand\mnfunc{u}
\newcommand\tmnfunc{w}

\newcommand\gfunc{g} 

\newcommand\crvmap{k}

\newcommand\HGMP{\mathcal{Y}(\manif,\Psp)}
\newcommand\HGMR{\mathcal{Y}(\manif,\RRR)}
\newcommand\HGMS{\mathcal{Y}(\manif,\Circle)}

\newcommand\smr{C^{\infty}(\manif,\RRR)}
\newcommand\sms{C^{\infty}(\manif,\Circle)}
\newcommand\smp{C^{\infty}(\manif,\Psp)}
\newcommand\smm{C^{\infty}(\manif,\manif)}

\newcommand\exlevel{L}

\newcommand\anbh{V}
\newcommand\bnbh{W}
\newcommand\hdif{h}

\newcommand\Gfld{G}

\newcommand\germm{C^{\infty}_{\pnt}(\manif)}

\newcommand\Jacidfz{\Delta(\mrsfunc,\pnt)}
\newcommand\Jacidfy{\Delta(\mrsfunc,y)}

\newcommand\sectMPP{\theta}
\newcommand\homood{\alpha}

\newcommand\pcr{e}
\newcommand\Pcr{E}
\newcommand\tpcr{\widetilde{\pcr}}
\newcommand\padj{q}
\newcommand\adjd{c}

\newcommand\Nbh{\mathcal{\nbh}}

\newcommand\dfval[2]{#1(#2)}
\newcommand\dfvalx[3]{\dfval{#1}{#2}(#3)}

\newcommand\condBndConst{{\rm (BC)}}
\newcommand\condFinCrVal{{\rm (FV)}}
\newcommand\condJacId{{\rm (J)}}

\newcommand\condE{{\rm (E)}}

\newcommand\pp{\pi}
\newcommand\arc{l}

\newcommand\prtwo{\prMPP}
\newcommand\prone{p_1}

\newcommand\pfact{\tau}

\newcommand\acnv{a}
\newcommand\lcnv{s}
\newcommand\dcnv{d}
\newcommand\kcnv{c}
\newcommand\Tn{T^{n}}

\newcommand\impcr{F_{n}(\Circle)}
\newcommand\hconf{\xi}

\newcommand\ConfSp{\FFF^{n}(\Circle)}
\newcommand\bij{\mu}

\newcommand\Gfunc{G}
\newcommand\factn{\nu}

\newcommand\Zk{Z_{\kcnv}}

\begin{document}
\begin{abstract}
Let $\manif$ be a smooth compact manifold and $\Psp$ be either $\RRR$ or $\Circle$.
There is a natural action of the groups $\DiffM$ and $\DiffM\times\DiffP$ on the space of smooth mappings $\smp$.
For $\mrsfunc\in\smp$ let $\StabfM$, $\StabfMP$, $\OrbfM$, and $\OrbfMP$ be the stabilizers and orbits of $\mrsfunc$ under these actions.
Recently, the author proved that under mild conditions on $\mrsfunc\in\smr$
the corresponding stabilizers and orbits are homotopy equivalent: $\StabfMR\sim\StabfM$ and $\OrbfMR\sim\OrbfM$.
These results are extended here to the actions on $\sms$.
It is proved that under the similar conditions (that are rather typical) we have that $\StabfMS\sim\StabfM$ and $\OrbfMS\sim\OrbfM\times\Circle$.
\end{abstract}
\maketitle

\newcommand\maksrefSectJcond{2}
\newcommand\maksrefSectML{5.1}
\newcommand\maksrefSectL{5.2}

\newcommand\maksrefPropDrDeRntwo{3.0.2}
\newcommand\maksrefLmOmOMR{6.0.1}

\newcommand\maksrefLmDrContr{3.0.1}
\newcommand\maksrefLmSmfunc{3.0.3}
\newcommand\maksrefLmKIsCont{4.0.1}
\newcommand\maksrefLmPrjDMRDR{5.0.1}
\newcommand\maksrefLmGluingVF{5.0.4}
\newcommand\maksrefLmFcompatVF{5.0.5}
\newcommand\maksrefLmSgmSmooth{5.2.1}
\newcommand\maksrefLmOmegaHom{5.2.2}

\section{Introduction}\label{sect:Intro}
In this note we extends the results of~\cite{Maks:StabR1} to the case of circle-valued mappings $\manif\to\Circle$.

Let $\manif$ be a smooth ($C^{\infty}$) connected compact $\dimM$-dimensional manifold, $\Psp$ be either the real line $\RRR$ or the circle $\Circle$,
and $\DiffM$ and $\DiffP$ be the groups of diffeomorphisms of $\manif$ and $\Psp$ respectively.

For $\pnt\in\manif$ let $\germm$ be the algebra of germs of smooth functions at $\pnt$.
If $\mrsfunc\in\germm$, then denote by $\Jacidfz$ the \emph{Jacobi ideal\/} of $\mrsfunc$ in $\germm$, i.e.\! the ideal generated by germs of partial derivatives of $\mrsfunc$ in some local coordinates at $\pnt$.

\begin{definition}
Let $\HGMP$ the subset of $\smp$ consisting of mappings $\mrsfunc$
satisfying the following conditions \condE\ and \condJacId:
\begin{enumerate}
\item[\condE]
$\mrsfunc$ is constant at every connected component of $\partial\manif$ and 
has only finitely many \emph{critical values};
\item[\condJacId]
for every critical point $y$ of $\mrsfunc$ there is a local representation $\mrsfunc:\RRR^{\dimM}\to\RRR$ in which the germ of the function $\mrsfunc(x)-\mrsfunc(y)$ belongs to the Jacobi ideal $\Jacidfy$.
\end{enumerate}
\end{definition}

\begin{remark}
The condition \condE\ is just an aggregate of the conditions \condBndConst\ and \condFinCrVal\ of \cite{Maks:StabR1}.
For the implications of the condition \condJacId\ see Section~\maksrefSectJcond\ of~\cite{Maks:StabR1}.
\end{remark}

If $\mrsfunc\in\smp$ satisfies \condE, then the set of critical points of $\mrsfunc$ may be infinite and there may be critical points on $\partial\manif$.
The values of $\mrsfunc$ on the connected components of $\partial\manif$ will be called \emph{boundary\/} ones.
All critical and boundary values of $\mrsfunc$ will be called \emph{exceptional} and their inverse images will be \emph{exceptional} levels of $\mrsfunc$.
Since $\manif$ is assumed to be compact, it follows that the set of exceptional values is finite.

Let us define the following groups.
If $\mrsfunc\in\smr$, then

$\bullet$~$\DiffRf\brf$ is the subgroup of $\DiffR$ consisting of diffeomorphisms that preserve orientation of $\RRR$, and leave the image of $\mrsfunc$ invariant;

$\bullet$~$\DiffRCrValf\brf$ is the subgroup of $\DiffRf\brf$ consisting of diffeomorphisms that also fix every exceptional value of $\mrsfunc$;

$\bullet$~$\DiffMR\brf = \DiffM \times \DiffRf\brf$.

\smallskip

If $\mrsfunc\in\sms$, then

$\bullet$~$\DiffSf$ is the group of preserving orientation diffeomorphisms of $\Circle$;

$\bullet$~$\DiffSCrValfInv\brf$ is the subgroup of $\DiffSf$ preserving the set of exceptional values of $\mrsfunc$;

$\bullet$~$\DiffSCrValf\brf$ is the (normal) subgroup of $\DiffSCrValfInv\brf$ fixing every exceptional value of $\mrsfunc$, thus
$\DiffSCrValfInv\brf/\DiffSCrValf\brf$ is a cyclic group $\ZZZ_{n}$, where $n$ is the number of exceptional values of $\mrsfunc$;

$\bullet$~$\DiffMS = \DiffM \times \DiffSf$.

\smallskip

For $\mrsfunc\in\smp$ we identify $\DiffM$ with the subgroup $\DiffM \times \id_{\Psp}$ of $\DiffMP$.
Then $\DiffM$ and $\DiffMP$ act on $\smp$ by the following formulas:
\chEq\begin{equation}\label{equ:act_M}
 \difM \cdot \mrsfunc  = \mrsfunc \circ \difM^{-1}
\end{equation}
\chEq\begin{equation}\label{equ:act_MP}
 (\difM,\difP) \cdot \mrsfunc  = \difP \circ \mrsfunc \circ \difM^{-1}.
\end{equation}
where $\mrsfunc\in\smp$, $\difM\in\DiffM$, and $\difP\in\DiffRf\brf$ or $\difP\in\DiffSf$ with respect to $\Psp$.
Let $\StabfM\brf$, $\StabfMP\brf$, $\OrbfM\brf$, and $\OrbfMP\brf$ be respectively the stabilizers and the orbits of $\mrsfunc\in\smp$ under the actions of $\DiffM$ and $\DiffMP$.
Evidently,
$\StabfM\brf\times\id_{\Psp}\subset\StabfMP\brf$ and $\OrbfM\brf\subset\OrbfMP\brf$.

We will often omit the notation $\mrsfunc$ and denote $\StabfM\brf$, $\OrbfM\brf$, $\DiffSCrValf\brf$ etc. simply by $\StabfM$, $\OrbfM$, $\DiffSCrValf$ and so on.
This will not lead to collisions.

Finally, we endow the spaces $\DiffM$, $\DiffP$, and $\smp$ with the corresponding $C^{\infty}$ Whitney topologies.
These topologies yield some topologies on $\DiffMP$ and on the
corresponding stabilizers and orbits of $\mrsfunc$.

\subsection{Main results}
Suppose that $\mrsfunc\in\HGMS$ has $n$ exceptional values and let $\prMPP:\DiffM\times\DiffS\to\DiffS$ be the standard projection.

\begin{theorem}\label{th:mainA}
$\DiffSCrValf\subset\prMPP(\StabfMS)$ and $\prMPP$ admits a section $\sectMPP:\DiffSCrValf\to\StabfMS$ being a homomorphism.
\end{theorem}
Thus putting $\pStabfMS=\prMPP^{-1}(\DiffSCrValf) \subset \StabfMS$, we get the following exact sequence split by $\sectMPP$: \
$ 1 \to \StabfM \to \pStabfMS \stackrel{\prMPP}{\to} \DiffSCrValf \to 1.$

\begin{corollary}
$\pStabfMS$ is a semi-direct product of $\StabfM$ and $\DiffSCrValf$.
\end{corollary}
\begin{corollary}
The embedding
$\StabfM \times\{\id_{\Circle}\} \subset \pStabfMS$ extends to a homeomorphism $\StabfM\times\DiffSCrValf\approx\pStabfMS$.
\end{corollary}

Notice that we have the following relations: 
\[
\StabfMS/\pStabfMS \ \stackrel{\prMPP}{\approx} \
\prMPP(\StabfMS)/\DiffSCrValf \ 
\subset \ \DiffSCrValfInv/\DiffSCrValf \ \approx \ZZZ_{n},
\]
whence 
\[
\prMPP(\StabfMS)/\DiffSCrValf  \approx \ZZZ_{\kcnv}
\]
for some $\kcnv$ that divides $n$.
Let $\dcnv$ be the index of $\prMPP(\StabfMS)$ in $\DiffSCrValfInv$, thus 
$n=\kcnv\dcnv$.

\begin{corollary}
$\StabfMS$ is homeomorphic with 
$\StabfM\times\DiffSCrValf \times \ZZZ_{\kcnv}$. 
\end{corollary}

\begin{definition}\label{defn:essent}
A critical point $\pnt$ of $\mrsfunc\in\smp$ is \emph{essential} if for  every neighborhood $\nbh$ of $\pnt$ there exists a neighborhood $\Nbh$ of $\mrsfunc$ in $\smp$ with $C^{\infty}$-topology such that every $\gfunc\in\Nbh$ has a critical point in $\nbh$.
\end{definition}

Let $\alpha:E\to\Circle$ be a $k$-dimensional vector bundle.
Then $E$ can be regarded as the quotient of $[0,1]\times\RRR^{k}$ by identifying $0\times\RRR^{k}$ with $1\times\RRR^{k}$ via some homeomorphism $\hdif:\RRR^{k}\to\RRR^{k}$.
The topological type of $E$ depends on the isotopy class of $\hdif$ only, whence we have two possibility for $E$.
If $\hdif$ preserves orientation, then $E\approx\Circle\times\RRR^{k}$.
Otherwise, we will say that $E$ is a ``twisted'' product and denote it by $\Circle\,\tilde\times\,\RRR^{k}$.

\begin{theorem}\label{th:mainB}
If $n=0$, then $\OrbfM=\OrbfMS$.
Otherwise, suppose that every critical level-set of $\mrsfunc$
includes either a connected component of $\partial\manif$ or an essential critical point or $\mrsfunc$.
Then the embedding $\OrbfM\subset\OrbfMS$ extends to a homeomorphism
\[
\begin{array}{ll}
\OrbfM \times  \Circle \,\tilde\times\,  
 \RRR^{n-1} \approx \OrbfMS, & \text{if \ $n$ \ is even and \ $\dcnv=n/\kcnv$ \ is odd}  \\
\OrbfM  \times 	 \Circle  \times \RRR^{n-1} \approx \OrbfMS,
& \text{otherwise}.
\end{array}
\]
\end{theorem}

\begin{remark}
Suppose that $n=0$. This means that $\mrsfunc:\manif\to\Circle$ has no critical points, whence this is a locally trivial fibering.
In this case we have:
$\prMPP(\StabfMS)=\DiffSf$,
the following sequence is exact:
\[ 1 \to \StabfM \to \StabfMS \stackrel{\prMPP}{\to} \DiffSf \to 1,\]
the projection $\prMPP$ admits a section $\sectMPP:\DiffSf\to\StabfMS$ being a homomorphism, $\StabfMS\approx\StabfM\times\DiffSf$ and $\OrbfM=\OrbfMS$.
\end{remark}

\begin{remark}
Suppose that $\DiffSCrValf\not=\prMPP(\StabfMS)$.
This means that there is $(\difM,\difR)\in\DiffMS$ such that $\difR\circ\mrsfunc=\mrsfunc\circ\difM$ and $\difR$ shifts exceptional values of $\mrsfunc$.
Then $\mrsfunc$ ``admits additional symmetry'' which seem not to be typical.
Thus for most mappings $\mrsfunc\in\HGMS$ we should have that $\DiffSCrValf=\prMPP(\StabfMS)$, whence 
\[\StabfMS\approx\StabfM\times\DiffSCrValf \quad  \text{and} \quad  
\OrbfMS\approx\OrbfM\times\Circle\times\RRR^{n-1}.\]
\end{remark}

The proofs of Theorems~\ref{th:mainA} and~\ref{th:mainB} will be given in sections~\ref{sect:proofA} and~\ref{sect:proofB}. 
They follow the line of~\cite{Maks:StabR1}.
Certain steps extends to $\Circle$-case almost literally,
therefore we often refer the reader to~\cite{Maks:StabR1}.

Notice that the results are similar to the $\RRR$-case in the part that 
the $\StabfM$ and $\OrbfM$ are the topological multiples of $\StabfMS$ and $\OrbfMS$ respectively.
The main cause of this is that $\Circle$ is a Lie group, whence there is a section of the evaluation map from $\DiffS$ to the $n$-th configuration space of $\Circle$, see Section~\ref{sect:DiffCr} and~\cite{FadellNeuwirth}.

\section{Configuration space $\ConfSp$ of $\Circle$}\label{sect:ConfSp}
We will regard $\Circle$ as the unit circle in the complex plane $\CCC$.
Recall that the $n$-th configuration space of $\Circle$ is the following subset 
\[
\ConfSp=\{ (x_0,\ldots,x_{n-1}) \in\Tn \ | \ x_i\not=x_j \ 
\text{for}\ i\not=j \}
\]
of $n$-dimensional torus $\Tn=\Circle\times\cdots\times\Circle$.

Let $\impcr$ be the connected component of $\ConfSp$ containing the point 
$(1, e^{2\pi i \frac{1}{n}}, \ldots , e^{2\pi i \frac{n-1}{n}}).$
Denote also 
\[
\symplex= \{ (x_{1},\ldots,x_{n-1}) \in \RRR^{n-1} \ | \
0 < x_1 < x_2  < \cdots < x_{n-1} < 1 \}.
\]
Then $\symplex$ is an open and convex subset of $\RRR^{n-1}$, whence it is diffeomorphic with $\RRR^{n-1}$.

\begin{lemma}\label{lm:Fn_S1_symplex}
$\impcr$ is diffeomorphic with $\Circle\times\symplex$.
\end{lemma}
\proof
Let $\pp:(0,1)\to\Circle$ be defined by the formula $\pp(t)=e^{2\pi i t}$. Evidently, $\pp$ 
is a diffeomorphism of $(0,1)$ onto $\Circle\setminus\{1\}$.
Consider now the following mapping $\hconf:\ConfSp\to\Circle\times\symplex$ defined by 
\[
\hconf(x_{0}, x_{1},\ldots, x_{n-1}) =
\bigl(
 x_{0}, \pp^{-1}(x_{1}/x_{0}), \ldots, \pp^{-1}(x_{n-1}/x_{0})
\bigr),
\]
where $x/y$ means the division of complex numbers.
Since $x_{\acnv}\not=x_{\acnv'}$ for $\acnv\not=\acnv'$, it follows that $x_{\acnv}/x_{0}\not=1$, whence $\pp^{-1}(x_{\acnv}/x_{0})$ is well-defined.
Notice also that $\hconf(\pnt_{0},\ldots,\pnt_{n-1})=(1,\frac{1}{n},\ldots,\frac{n-1}{n})$.
It is easy to see that $\hconf$ is a surjective local diffeomorphism. Moreover, it admits a section $s:\Circle\times\symplex\to\ConfSp$ defined by the formula
\[s(\phi, t_1,\ldots,t_{n-1})=(\phi, \phi\cdot\pp(t_1),\ldots,\phi\cdot\pp(t_{n-1}))\]
such that $s(1,\frac{1}{n},\ldots,\frac{n-1}{n})=(\pnt_{0},\ldots,\pnt_{n-1})$.
Hence $\hconf$ yields a diffeomorphism between the connected component $\impcr$ of $(\pnt_{0},\ldots,\pnt_{n-1})$ in $\ConfSp$ and $\Circle\times\symplex$.
\endproof
Consider the action of the group $\ZZZ_{n}$ on $\ConfSp$ by cyclic permutation 
$\sigma:\ConfSp\to\ConfSp$ of coordinates:
\[\sigma(x_{0},\ldots,x_{n-1})=(x_{1},\ldots,x_{n-1},x_{0}).\]
Evidently, this action is free and $\sigma(\impcr)=\impcr$.
Moreover, $\sigma$ preserves orientation of $\impcr$ iff and only if $n$ is odd.

Suppose that $n=\kcnv\dcnv$ and let $\Zk$ be a (cyclic) subgroup (of order $\kcnv$) of $\ZZZ_{n}$ generated by $\sigma^{\dcnv}$.

\begin{lemma}\label{lm:FnZk_S1_symplex}
Suppose that $n$ is even and $\dcnv$ is odd, then $\impcr/\Zk$ is diffeomorphic with a ``twisted'' product $\Circle\,\tilde\times\,\symplex$.
Otherwise, we have $\impcr/\Zk\approx\Circle\times\symplex$.
\end{lemma}
\proof
The proof is direct and based on the remark that 
$\sigma^{\dcnv}$ reverses orientation if and only if $n$ is even and $\dcnv$ is odd.
\endproof

\section{The groups $\DiffSf$ and $\DiffSCrValf$}\label{sect:DiffCr}
Suppose that $\mrsfunc\in\HGMS$ and 
$\pnt_{\acnv}=e^{2\pi i \frac{\acnv}{n}}$ $(\acnv=0,\ldots,n-1)$ are all of the exceptional values of $\mrsfunc$. We will assume that $n\geq1$.
Then $\DiffSCrValf$ is the subgroup of $\DiffSf$ that fixes every point $\pnt_{\acnv}$.

Let $\pcr:\DiffSf\to\ConfSp$ be the \emph{evaluation mapping\/} defined by
\chEq\begin{equation}\label{equ:eval_map}
\pcr(\difR)=\left( \,
\difR(\pnt_{0}), \ldots, \difR(\pnt_{n-1}) \right),
\quad \difR\in\DiffSf.
\end{equation}
Then the image $\pcr(\DiffSf)$ coincides with the connected component $\impcr$ of $\ConfSp$ containing $(\pnt_{0},\ldots,\pnt_{n-1})$.

Notice that $\pcr$ is constant on the adjacent classes of $\DiffSf$ by $\DiffSCrValf$ and yields a bijective continuous mapping $\tpcr:\DiffSf/\DiffSCrValf\to\impcr$ such that
\[
\begin{CD}
\pcr=\tpcr\circ\padj:\DiffSf @>{\padj}>>\DiffSf/\DiffSCrValf @>{\tpcr}>> \impcr,
\end{CD}
\]
where $\padj$ is a factor-mapping.

By Lemma~\ref{lm:Fn_S1_symplex} we have a diffeomorphism $\hconf:\impcr\approx\Circle\times\symplex$.
Thus in order to show that $\tpcr$ is a homeomorphism it suffices to prove the following statement.

\begin{proposition}\label{pr:DS_DSF_S1_symplex}
The mapping $\hconf\circ\tpcr:\DiffSf/\DiffSCrValf\to\Circle\times\symplex$ is a homeomorphism in all Whitney topologies.
\end{proposition}
\begin{proof}
Consider the composition:
\[
 \Pcr =\hconf\circ\tpcr\circ\padj:
 \DiffS \ \stackrel{\padj}{\to} \
 \DiffSf/\DiffSCrValf \ \stackrel{\tpcr}{\to} \
 \impcr \ \stackrel{\hconf}{\to} \
 \Circle\times\symplex.
\]
By the arguments used in Proposition~\maksrefPropDrDeRntwo\ of~\cite{Maks:StabR1} it suffices to construct a continuous section $\sct:\Circle\times\symplex\to\DiffSf$ of $\Pcr$.
The following statement is Lemma~\maksrefLmSmfunc\ of~\cite{Maks:StabR1} with a slight adaptation of notations to our case:
\begin{lemma}\label{clm:constr_section_not_onto}
There exists a smooth function $\mnfunc:\RRR\times\symplex\to\RRR$
with the following properties:
\begin{enumerate}
\em\item
$\mnfunc'_{t}(t; x_1,\ldots,x_{n-1})>0$ for all $(t;x_1,\ldots,x_{n-1}) \in\RRR\times\symplex$;
\item
$\mnfunc(t;x_1,\ldots,x_{n-1}) = t$ for $t\leq 0$ and $t\geq 1$;
\item
$\mnfunc\left(\frac{k}{n}, x_1,\ldots,x_{n-1} \right) = x_k$ for $k=1,2,\ldots,n-1$.
\item
$\mnfunc\left( t; \frac{1}{n},\frac{2}{n},\ldots,\frac{n-1}{n}\right) = t$ for $t\in\RRR$.
\end{enumerate}
\end{lemma}

Let $\mnfunc$ be a function of this lemma.
Then $\mnfunc$ maps the set $[0,1]\times\symplex$ onto $[0,1]$ so that
$\mnfunc(0, x)=0$ and $\mnfunc(1, x)=1$.
Hence $\mnfunc$ yields a \emph{continuous\/} mapping
$\tmnfunc:\Circle\times\symplex\to\Circle$ defined by factorization of $[0,1]$ to $\Circle$ by the mapping $\pp$.

Moreover,
$\mnfunc'_{t}(0, x)=\mnfunc'_{t}(1, x)=1$ and
$\mnfunc^{(s)}_{t}(0, x)=\mnfunc^{(s)}_{t}(1, x)=0$ for $s\geq 2$ and $x\in\symplex$.
Hence, $\tmnfunc$ is in fact $C^{\infty}$.

Similarly to Lemma~\maksrefLmSmfunc\ of~\cite{Maks:StabR1} it can be proved that the following mapping $\sct:\Circle\times\symplex\to\DiffSf$ defined by:
\chEq\begin{equation}\label{equ:sect}
\sct(x_1,\ldots,x_{n-1})(t) = \mnfunc(t;x_1,\ldots,x_{n-1})
\end{equation}
is a section of $\Pcr$.
This completes our proposition.
\end{proof}

\section{Proof of Theorem~\ref{th:mainA}}\label{sect:proofA}
Suppose that $\mrsfunc\in\sms$ satisfies the condition \condE.
Let $n$ be the number of exceptional levels of $\mrsfunc$.

Notice that for $n=0$ we have $\DiffSCrValf=\DiffSCrValfInv=\DiffSf$.

The following statement are analogues of Lemmas~\maksrefLmDrContr\ and~\maksrefLmPrjDMRDR\ of~\cite{Maks:StabR1}.
We leave the proof to the reader.
\begin{lemma}\label{lm:DiffSCrValf_is_contr}
If $n\geq1$ then the group $\DiffSCrValf$ is contractible. \qed
\end{lemma}

\begin{lemma}\label{lm:pStMP_in_Pcrf}
Let $\prMPP:\DiffMS\to\DiffSf$ be the projection onto the second multiple.
Then $\prMPP(\StabfMS)\subset\DiffSCrValfInv$. \qed
\end{lemma}

Suppose that $\mrsfunc\in\HGMS$.
We have to construct a continuous section $\sectMPP:\DiffSCrValf\to\StabfMS$ of $\prMPP$ being a homomorphism.
The proof follows the line of the $\RRR$-case, but there are some exceptions. 
Let us briefly recall the idea.

Notice that every section $\sectMPP:\DiffSCrValf\to\StabfMS$ of $p$ should be of the form $\sectMPP(\difR) = (\difShift(\difR),\difR)$,
where $\difShift:\DiffSCrValf\to\DiffM$ is a continuous homomorphism such that
$\difR\circ\mrsfunc=\mrsfunc\circ\difShift(\difR)$.
Such a section was constructed in~\cite{Maks:StabR1} via a vector field of a special form.

\begin{definition}\label{defn:f-compat}
Let $\mrsfunc\in\smp$, where $\Psp$ is either $\RRR$ of $\Circle$.
We say that a vector field $\Fld$ on $\manif$ is \emph{$\mrsfunc$-compatible\/} if the following conditions (1) and (2) hold true:

(1) if $y\in\manif$ belongs to a non-exceptional level-set of $\mrsfunc$, then
\[\dfvalx{d\mrsfunc}{\Fld}{y}\not=0;\]

(2) if $y$ belongs to an exceptional level-set $\exlevel=\mrsfunc^{-1}(\mrsfunc(y))$, then
\[\dfvalx{d\mrsfunc}{\Fld}{x}=\eps(\mrsfunc(x)-\mrsfunc(y)),\]
for all $x$ from some neighborhood of $y$ and some $\eps=\pm1$.
Evidently $\eps$ is the same for all points of $\exlevel$.
Therefore we call $\exlevel$ \emph{attractive} if $\eps=-1$ and \emph{reflective} otherwise.
\end{definition}

Let $\mrsfunc\in\HGMR$. We proved in Lemma~\maksrefLmFcompatVF\ of~\cite{Maks:StabR1} that there exists an $\mrsfunc$-compatible vector field $\Fld$.
Let $\flow$ be the flow generated by $\Fld$.
Then we defined the {\em shift\/} mapping $\Shift:\smr \to \smm$ along trajectories of $\flow$ by the following formula, see~\cite{Maks:Shifts}:
\[ \Shift(\afunc)(x)=\flow(x,\afunc(x)),\qquad x\in\manif, \ \afunc\in\smr, \]
and the mapping $\claimbmap:\DiffRCrValf\to\smr$ such that the composition
\[
\begin{CD}
\difShift = \Shift \circ \claimbmap:
  \DiffRCrValf @>{\claimbmap}>> \smr @>{\Shift}>> \smm
\end{CD}
\]
is what we need, i.e.\! a continuous homomorphism such that for every $\difP\in\DiffRCrValf$ we have
$(\difShift(\difP), \difP) \in \StabfMR$:
\[
\difP\circ\mrsfunc(x) \ = \ \mrsfunc\circ \flow(x,\claimbmap(\difP)(x)) \
= \ \mrsfunc\circ \difShift(\difP)(x), \qquad x\in\manif.
\]

\begin{lemma}\label{lm:f-comp_notexists}
Suppose that $\mrsfunc\in\sms$ satisfies the condition \condE.
Then $\mrsfunc$-compatible vector field may exist only if the number of exceptional levels of $\mrsfunc$ is \emph{even}.
\end{lemma}
\proof
If $\Fld$ is an $\mrsfunc$-compatible vector field then it follows from the definition that attractive and reflective levels must alternate with each other.
\endproof

\begin{example}
Suppose that $\manif=\Circle$ and let $\mrsfunc:\manif\to\Circle$ be a smooth homeomorphism with $4$ critical points $\pnt_0,\pnt_1,\pnt_2,\pnt_3$, 
see Figure~\ref{fig:compat_vf}, where $\mrsfunc$ is an ``almost central projection'' of the external circle to the internal one.
Suppose that $\Fld$ is an $\mrsfunc$-compatible vector field on $\Circle$. Then the trajectories of $\Fld$ should be of the form shown in the right part of Figure~\ref{fig:compat_vf}.
\end{example}

Let $\mrsfunc\in\HGMS$. 
If $n$ is odd, then $\Fld$ and therefore $\flow$, $\Shift$, and $\claimbmap$ \emph{may exist only locally}.
Nevertheless, we will show that it is possible to define them so that the composition $\difShift = \Shift \circ \claimbmap$ is well-defined on all of $\manif$.

\subsubsection{Case $n=0$.}
In this case, following Lemma~\maksrefLmFcompatVF\ of~\cite{Maks:StabR1}, we can define a \emph{global\/} $\mrsfunc$-compatible vector field $\Fld$ on $\manif$ such that $d\mrsfunc(\Fld)>0$.
Then the construction of $\difShift$ is similar to~\cite{Maks:StabR1}. 

\subsubsection{Case $n\geq1$.}
We can assume that $\pnt_{\acnv} = e^{2\pi i \frac{\acnv}{n}}$ $(\acnv=0,\ldots,n-1)$ are all of the exceptional values of $\mrsfunc$.
Below all indexes will be taken modulo $n$.
Let $\exlevel_{\acnv}=\mrsfunc^{-1}(\pnt_{\acnv})$ and
$\exlevel=\mathop\cup\limits_{\acnv=0}^{n-1}\exlevel_{\acnv}$.

If $n=1$ then put $\arc_{0}=\Circle\setminus\pnt_{0}$.
Otherwise, $n\geq2$, therefore the points
$\pnt_{\acnv}$ and $\pnt_{\acnv+1}$ divides $\Circle$ into two open arcs.
Let $\arc_{\acnv}$ be that arc which starts from $\pnt_{\acnv}$ and finishes at $\pnt_{\acnv+1}$ with respect to the positive orientation of $\Circle$.
Denote also $\bnbh_{\acnv} = \mrsfunc^{-1}(\arc_{\acnv})$.

Following Lemma~\maksrefLmFcompatVF\ of~\cite{Maks:StabR1} for every $\acnv=0,\ldots,n-1$ we can define
\begin{enumerate}
\item
an $\mrsfunc$-compatible vector field $\Gfld_{\acnv}$ on $\bnbh_{\acnv}$
such that $d\mrsfunc(\Gfld(x))>0$ for $x\in\bnbh_{\acnv}$;
\item
a neighborhood $\nbh_{\acnv}$ of $\exlevel_{\acnv}$ and an $\mrsfunc$-compatible vector field $\Fld_{\acnv}$ on $\nbh_{\acnv}$
such that $\nbh_{\acnv}\cap\nbh_{\acnv'}=\varnothing$ for $\acnv\not=\acnv'$
and 
$\mrsfunc(x)-\mrsfunc(y)=d\mrsfunc(\Fld_{\acnv})(x)$ for $x\in\nbh_{\acnv}$ and $y\in\exlevel_{\acnv}$.
\end{enumerate}

Let $\nbh_{\acnv}^{-}$ and $\nbh_{\acnv}^{+}$ be the ``lower'' and ``upper'' parts of $\nbh_{\acnv}$ with respect to the orientation of $\Circle$, i.e.
\[
\nbh_{\acnv}^{\pm} = \{ x \in \nbh_{\acnv} \ | \ 
\pm d\mrsfunc(\Fld_{\acnv}(x))>0
\}
\]
Gluing $\Fld_{\acnv}$ with $\Gfld_{\acnv}$ on $\nbh_{\acnv}^{+}$ we can assume that
$\Fld_{\acnv}=\Gfld_{\acnv}$ on $\nbh_{\acnv}^{+}$.
Similarly, gluing $\Fld_{\acnv}$ with $-\Gfld_{\acnv-1}$ on $\nbh_{\acnv}^{-}$
we can also assume that $\Fld_{\acnv}=-\Gfld_{\acnv-1}$ on $\nbh_{\acnv}^{-}$, see Lemma~\maksrefLmGluingVF\ of~\cite{Maks:StabR1}.

 \chFig\begin{figure}[ht]
 \includegraphics[height=5cm]{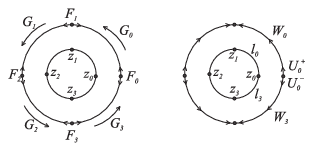}
 \caption{}\protect\label{fig:compat_vf}
 \end{figure}

Let $\aflow_{\acnv}:\bnbh_{\acnv}\times\RRR\to\bnbh_{\acnv}$ be the flow that generated by $\Gfld_{\acnv}$.
Then, see~\S\maksrefSectML\ of~\cite{Maks:StabR1}, for every $\difP\in\DiffSCrValf$ there is a smooth function $\clmdif'_{\difP}:\bnbh_{\acnv}\to\RRR$ such that
\[
\difR\circ\mrsfunc(x)=
\mrsfunc\circ\aflow_{\acnv}(x,\clmdif_{\difP}(x)).
\]

Let also $\flow_{\acnv}:\anbh_{\acnv}\times(-\eps,\eps)\to\nbh_{\acnv}$ be the local flow generated by $\Fld_{\acnv}$ on some neighborhood $\anbh_{\acnv}\subset\nbh_{\acnv}$ of $\exlevel_{\acnv}$.
Then, see \S\maksrefSectL\ of~\cite{Maks:StabR1}, for every $\difP\in\DiffSCrValf$ there is a smooth function $\clmdif_{\difP}$ defined near $\exlevel_{\acnv}$ in $\anbh_{\acnv}$ such that
\[
\difR\circ\mrsfunc(x)=
\mrsfunc\circ\flow_{\acnv}(x,(x)).
\]
These flows and functions are compatible in the following sense:
let $\difR\in\DiffSCrValf$ and suppose that $\clmdif_{\difP}$ is defined at
$x\in\nbh_{\acnv}$. Then

\begin{enumerate}
\item[\rm(a)]
if $x\in\nbh_{\acnv}^{+}$, then $\Fld_{\acnv}(x)=\Gfld_{\acnv}(x)$, whence $\flow_{\acnv}(x,t)=\aflow_{\acnv}(x,t)$ and $\clmdif_{\difP}(x)=\clmdif'_{\difP}(x)$.
\item[\rm(b)]
if $x\in\nbh_{\acnv}^{-}$, then $\Fld_{\acnv}(x)=(-1)^{n}\,\Gfld_{\acnv}(x)$, therefore we have that 
$\flow_{\acnv}(x,t)=\aflow_{\acnv}(x,(-1)^{n}\,t)$ and $\clmdif_{\difP}(x)=(-1)^{n}\,\clmdif'_{\difP}(x)$.
\end{enumerate}
In both cases we see that 
$\flow_{\acnv}(x,\clmdif_{\difP}(x)) = \aflow_{\acnv}(x,\clmdif'_{\difP}(x))$.
Then the following mapping $\difShift(\difP):\manif\to\manif$ given by the formula
\[
\difShift(\difP)(x) = \left\{
\begin{array}{ll}
\flow_{\acnv}(x,\clmdif_{\difP}(x)), & x\in\nbh_{\acnv} \\
\aflow_{\acnv}(x,\clmdif'_{\difP}(x)), & x\in\anbh_{\acnv}
\end{array}
\right.
\]
is well-defined on all of $\manif$ and satisfies the condition
$\difR\circ\mrsfunc = \mrsfunc\circ\difShift(\difP)$.

Similarly to Lemmas~\maksrefLmSgmSmooth\ and~\maksrefLmOmegaHom\ of~\cite{Maks:StabR1} it can be shown that $\difShift(\difP)$ is a diffeomorphism of $\manif$ and
the mapping $\difShift:\DiffSCrValf\to \DiffM$ defined by $\difR\mapsto\difShift(\difR)$ is a continuous homomorphism.
Hence the correspondence
$\difR \mapsto (\difShift(\difR),\difR)$ is a continuous section $\sectMPP:\DiffSCrValf\to\StabfMS$ of $\prtwo$.

Theorem~\ref{th:mainA} is completed.

\begin{remark}
If $n\geq2$ is \emph{even}, then 
we could replace $\Fld_{\acnv}$ with $-\Fld_{\acnv}$ and $\Gfld_{\acnv}$ with $-\Gfld_{\acnv}$ for even $\acnv=0,2,\ldots,n-2$, see Figure~\ref{fig:compat_vf}.
This would give us a global $\mrsfunc$-compatible vector field $\Fld$ and therefore we could define $\difShift$ similarly to~\cite{Maks:StabR1}.
But we preferred not to separate the cases of even and odd $n$.
\end{remark}

\begin{corollary}\label{cor:exists_difR_periodic}
Suppose that $(\difM,\difR)\in\StabfMS$ but $\difR\not\in\DiffSCrValf$.
Then there is a path in $\StabfMS$ between $(\difM,\difR)$ and $(\bar\difM,\bar\difR)$ such that $\bar\difR$ is periodic and $\bar\difR^{-1}\circ\difR\in\DiffSCrValf$.
\end{corollary}
\proof
We have 
$\difR\circ\mrsfunc=\mrsfunc\circ\difM$.
Since $\difR\not\in\DiffSCrValf$, it follows that $\difR$ cyclically shifts $\pnt_{\acnv}$, thus $\difR(\pnt_{\acnv}) = \pnt_{\acnv+\lcnv}$ for some $\lcnv\in\ZZZ_{\kcnv}$ and every $\acnv=0,\ldots,n-1$, where indexes are taken modulo $n$.

Let $\bar\difR(z)=e^{2\pi i \lcnv/n}\cdot z$ be the rotation of $\Circle$ by the angle $2\pi\lcnv/n$.
Then $\bar\difR\circ\difR^{-1}\in\DiffSCrValf$, whence
$\bar\difR\circ\difR^{-1}\circ \mrsfunc = \mrsfunc \circ  \difShift(\bar\difR\circ\difR^{-1})$.
Therefore
\[
\bar\difR\circ\mrsfunc =
(\bar\difR\circ\difR^{-1})\circ\difR\circ\mrsfunc=
(\bar\difR\circ\difR^{-1})\circ\mrsfunc\circ\difM=
\mrsfunc\circ\difShift(\bar\difR\circ\difR^{-1})\circ\difM.
\]
Thus putting $\bar\difM=\difShift(\bar\difR\circ\difR^{-1})\circ\difM$
we see that $(\bar\difM,\bar\difR)\in\StabfMS$ and $\bar\difR$ is periodical.

By Lemma~\ref{lm:DiffSCrValf_is_contr} $\DiffSCrValf$ is contractible. Therefore there exists an isotopy 
$\difR_{t}:\Circle\to\Circle$ between $\difR_{0}=\bar\difR\circ\difR^{-1}$ and $\difR_{1}=\id_{\Circle}$.
It gives a path $\omega_{t}=(\difShift(\difR_{t})\circ\difM,\difR_{t}\circ\difR)$
between 
$\omega_{0} = (\difShift(\bar\difR\circ\difR^{-1})\circ\difM,\bar\difR)=
(\bar\difM,\bar\difR)$ and $\omega_{1} = (\difM,\difR)$.
\endproof

\section{Proof of Theorem~\ref{th:mainB}}\label{sect:proofB}

Suppose that $n=0$. In this case $\mrsfunc$ is a locally trivial fibering over $\Circle$. Then $\DiffSCrValf=\DiffSf$, therefore $\difShift$ is defined on all of $\DiffSf$.
We have to show that $\OrbfM=\OrbfMS$.
Evidently, $\OrbfM\subset\OrbfMS$.
Conversely, if $\gfunc=\difR\circ\mrsfunc\circ\difM^{-1}\in\OrbfMS$ for $\difR\in\DiffSf$ and $\difM\in\DiffM$, then we also have $\gfunc=\mrsfunc\circ\difShift(\difR)\circ\difM^{-1}\in\OrbfM$.
Thus $\OrbfM=\OrbfMS$.

Now let $n\geq 1$. We will assume that $\pnt_{\acnv}=e^{2\pi i \frac{\acnv}{n}}$ $(\acnv=0,\ldots,n-1)$ are all of the exceptional values of $\mrsfunc$.
Then we have a well-defined mapping \[\crvmap:\OrbfMS\to\impcr/\ZZZ_{n}\] corresponding to every $\gfunc\in\OrbfMS$ the cyclically ordered set of its exceptional values.
The following lemma is an analogue of Lemma~\maksrefLmKIsCont\ of~\cite{Maks:StabR1}.
The proof is quite the same and we leave it to the reader.
\begin{lemma}\label{lm:cond_crlev_cont}
Suppose that every critical level-set of $\mrsfunc$
includes either a connected component of $\partial\manif$ or an essential critical point or $\mrsfunc$.
Then the mapping $\crvmap$ is continuous in $C^{\infty}$-topology of $\OrbfMS$.\qed
\end{lemma}

Consider the mapping $\pfact:\DiffMS\to\OrbfMS\times\impcr$ defined by:
\[
\pfact(\difM,\difR) = \bigl( \difR\circ\mrsfunc\circ\difM^{-1}; 
\difR(\pnt_{0}),\ldots,\difR(\pnt_{n-1}) \bigr).
\]
Evidently, $\pfact$ is continuous.
Denote its image by $\pOrbfMS$: \[\pOrbfMS:=\pfact(\DiffMS) \subset\OrbfMS\times\impcr.\]

Recall that we have the action of $\ZZZ_{n}$ on $\impcr$ generated by the cyclic permutation $\sigma$ of coordinates.
This action together with the trivial action on $\OrbfMS$ yields an action of $\ZZZ_{n}$ on $\OrbfMS\times\impcr$.

By Theorem~\ref{th:mainA}, $\prMPP(\StabfMS)/\DiffSCrValf$ is a cyclic group of some order $\kcnv$ dividing $n$.
Denote $\dcnv=n/\kcnv$. Then we can identify it with the subgroup $\Zk$ of $\ZZZ_{n}$ generated by $\sigma^{\dcnv}$.

Let also $\prone:\pOrbfMS\to\OrbfMS$ the restrictions of standard projection $\OrbfMS\times\impcr\to\OrbfMS$ to $\pOrbfMS$.
\begin{lemma}\label{lm:struct_pOrbfMS}
$\pOrbfMS$ is invariant under $\Zk$ and 
$\prone$ 
yields a continuous bijection $\bij:\pOrbfMS/\Zk\to\OrbfMS$.
If $\crvmap$ is continuous, then $p_1$ is a $\kcnv$-sheet covering and $\bij$ is a homeomorphism.
\end{lemma}
\proof
In order to simplify notations we denote a point $(x_{0},\ldots,x_{n-1})\in\impcr$ by $\{x_{\acnv}\}$.

Let $(\gfunc, \{x_{\acnv}\})\in\pOrbfMS$, i.e.
$\gfunc=\difR\circ\mrsfunc\circ\difM^{-1}$ and $x_{\acnv}=\difR(\pnt_{\acnv})$.
We have to prove that 
\[
\lcnv\cdot(\gfunc, \{x_{\acnv}\})=
(\gfunc,\{x_{\dcnv\lcnv + \acnv}\})\in\pOrbfMS
\]
 for every $\lcnv\in\Zk$, i.e.\! $\gfunc=\hat\difR\circ\mrsfunc\circ\hat\difM^{-1}$ for some $(\hat\difM,\hat\difR)\in\DiffMS$ such that 
$\hat\difR(\pnt_{\acnv}) = x_{\dcnv\lcnv + \acnv}$.
Let $\bar\difR\in\DiffSf$ be the rotation of $\Circle$ by the angle $2\pi/\kcnv$.
Then $\bar\difR^{\lcnv}(\pnt_{\acnv}) = \pnt_{\dcnv\lcnv + \acnv}$.
Using Corollary~\ref{cor:exists_difR_periodic}, for every $\lcnv\in\Zk$ choose $\difM_{\lcnv}\in\DiffM$ such that 
$(\difM_{\lcnv}, \bar\difR^{\lcnv})\in \StabfMS$, i.e.
$\mrsfunc =  \bar\difR^{\lcnv}\circ \mrsfunc\circ \difM_{\lcnv}^{-1}$.
Then
\[
\gfunc=\difR\circ\mrsfunc\circ\difM^{-1} = 
\difR\circ\bar\difR^{\lcnv}\circ \mrsfunc\circ \difM_{\lcnv}^{-1}\circ\difM^{-1}.
\]
Moreover, 
$\difR\circ\bar\difR^{\lcnv}(\pnt_{\acnv}) =
\difR(\pnt_{\dcnv\lcnv + \acnv})=x_{\dcnv\lcnv + \acnv}$.
Thus putting $\hat\difR=\difR\circ\bar\difR^{\lcnv}$ and 
$\hat\difM=\difM\circ\difM_{\lcnv}$ we obtain that 
\[
(\gfunc,\{x_{\dcnv\lcnv + \acnv}\})=\pfact(\hat\difM,\hat\difR)\in \,\text{image of $\pfact$}\, \equiv \, \pOrbfMS.
\]
Thus $\pOrbfMS$ is invariant under $\Zk$.
Since $\Zk$ is finite and its action is free and $p_1$-equivariant, it follows that the factor mapping $\pOrbfMS\to\pOrbfMS/\Zk$ is a $\kcnv$-sheet covering.
Moreover, since $p_1$ is onto, we obtain that $p_1$ yields a continuous bijection $\bij:\pOrbfMS/\Zk\to\OrbfMS$.

Suppose that $\crvmap$ is continuous.
We will show that $\prone$ is a local homeomorphism. This will imply that so is $\bij$, whence $\bij$ is in fact a homeomorphism.
It suffices to prove that $\prone$ admits continuous local sections, i.e.\!
for every $(\gfunc,x)\in\pOrbfMS$ there exists a neighborhood $\Nbh_{\gfunc}$ of $\gfunc$ in $\OrbfMS$ and a continuous mapping $\Gfunc:\Nbh_{\gfunc}\to\impcr$ such that $(\gfunc',\Gfunc(\gfunc'))\in\pOrbfMS$ and $\Gfunc(\gfunc)=x$.

Notice that we have the following maps:
\[
\impcr \stackrel{\factn}{\longrightarrow} \impcr/\ZZZ_{n} \stackrel{\crvmap}{\longleftarrow} \OrbfMS.
\]

Let $(\gfunc,x)\in\pOrbfMS$ and $[x]=\factn(x)$ be the corresponding class of $x$ in $\impcr/\ZZZ_{n}$. 
Then $\crvmap(\gfunc)=[x]$.
Since $\factn$ is covering, there is a neighborhood $\nbh_{x}$ of $x$ in $\impcr$ and a neighborhood $\anbh_{[x]}$ of $[x]$ in $\impcr/\ZZZ_{n}$ such that $\factn$ homeomorphically maps $\nbh_{x}$ onto $\anbh_{[x]}$.
Let $\Nbh_{\gfunc}=\crvmap^{-1}(\anbh_{[x]})$. 
Since $\crvmap$ is continuous, we obtain that $\Nbh_{\gfunc}$ is an open neighborhood of $\gfunc$. Then the mapping 
$\factn^{-1}\circ\crvmap:\Nbh_{\gfunc}\to\nbh_{x}$ is a local section of $\prone$.
\endproof

Consider now the projection \[p_2:\pOrbfMS\subset\OrbfMS\times\impcr\to\impcr.\]
Then the composition 
\[
\begin{CD}
p_2\circ\pfact:\DiffMS @>{\pfact}>> \OrbfMS\times\impcr @>{p_2}>>\impcr.
\end{CD}
\]
is given by $(\difM,\difR)\mapsto\bigl(\difR(\pnt_0),\ldots,\difR(\pnt_{n-1})\bigr)$.
Thus we have the following commutative diagram:
\[
\begin{CD}
\StabfM       @>{}>> \DiffM   @>{}>>      \DiffM/\StabfM        @>{}>>      \OrbfM  \\
 @V{}VV              @V{}VV                @V{}VV                            @V{}VV \\
\pStabfMS      @>{}>> \DiffMS  @>{}>>      \DiffMS/\pStabfMS      @>{}>>      \pOrbfMS \\ 
@V{\prMPP}VV         @V{}VV                @V{}VV                    @V{p_2}VV \\
\DiffSCrValf  @>{}>> \DiffSf  @>{\padj}>> \DiffSf/\DiffSCrValf  @>{\adjd}>> \impcr
\end{CD}
\]
in which $\adjd$ is a homeomorphism and the other right horizontal arrows are continuous bijections, upper vertical and left horizontal arrows are embeddings, and lower vertical ones are surjective mappings.
Notice that 
\begin{enumerate}
\item
$\adjd\circ\padj$ admit a continuous section $\sct$ being a homomorphism (Proposition~\ref{pr:DS_DSF_S1_symplex}), 
\item
$\prtwo$ admit a continuous section $\sectMPP$ (Theorem~\ref{th:mainA}), and
\item $p_2$ is continuous.
\end{enumerate}
Then it follows by the arguments of Proposition~\maksrefLmOmOMR\ of~\cite{Maks:StabR1}, that the embedding $\OrbfM \equiv \OrbfM \times (\pnt_{0},\ldots,\pnt_{n-1}) \subset\pOrbfMS$ extends to a homeomorphism
\[
\homood:\OrbfM\times\impcr \approx \pOrbfMS, \qquad
\homood(\gfunc,x)=(\sct(x)\circ\gfunc,x),
\]
where $\gfunc\in\OrbfM$ and $x\in\impcr$.

Finally, by Lemma~\ref{lm:struct_pOrbfMS} we have that $\OrbfMS=\pOrbfMS/\Zk$.
Since $\Zk$ trivially acts on $\OrbfM$, we obtain a homeomorphism:
\[
\OrbfMS \approx \OrbfM\times(\impcr/\Zk).
\]
Then it remains to apply Lemma~\ref{lm:FnZk_S1_symplex}.
\qed


\end{document}